\documentclass[twoside,english]{elsarticle}
\usepackage[LGR,T1]{fontenc}
\usepackage[latin9]{inputenc}
\usepackage{array}
\usepackage{verbatim}
\usepackage{rotating}
\usepackage{float}
\usepackage{calc}
\usepackage{units}
\usepackage{bbding}
\usepackage{multirow}
\usepackage{amsthm}
\usepackage{amsbsy}
\usepackage{amstext}
\usepackage{amssymb}
\usepackage{arcs}
\usepackage{graphicx}
\usepackage{xcolor}
\usepackage{caption}
\usepackage{soul}
\setul{1pt}{0.8pt}
\makeatletter

\DeclareFontEncoding{LGR}{}{}
\DeclareTextSymbol{\~}{LGR}{126}
\providecommand{\tabularnewline}{\\}
\theoremstyle{plain}

\theoremstyle{plain}

\journal{Example: Nuclear Physics B}

\makeatother

\usepackage{babel}
\providecommand{\conjecturename}{Conjecture}
\providecommand{\theoremname}{Theorem}

\begin{document}
\journal{{\tt arXiv}\it}

\title{Parafermi algebra and semiprimality\tnoteref{t1,t2}}

\tnotetext[t1]{This document builds on aspects addressed in an earlier study titled
``Parafermi algebra and interordinality''\,\,\cite{Merkel}.}

\tnotetext[t2]{In ``Parafermi algebra and interordinality'', the central theme
was the implications raised by the special case that two parafermi
algebras are of Mersenne-wise neighboring orders. The present document
is meant to be largely self-contained, but an in-depth study of the
earlier work is helpful and therefore recommended.}

\author{U.~Merkel\corref{cor1}}

\ead{merkel.u8@googlemail.com}

\cortext[cor1]{Corresponding author}

\address[rvt]{Universit{\"a}tsstr. 38, 70569 Stuttgart, Germany}

\address{}
\begin{abstract}
The principles of  rational parafermi algebra of Mersenne order $p=2^n-1$, developed in an earlier work \cite{Merkel}, are recapitulated in the introduction:  A Continued  Kronecker Square-root process (CKS)  dispenses with irrational numbers.The process creates two types of lower triangular square matrices, called of $G$-type and $J$-type, whose elements  decompose Catalan numbers in particular ways.  Trilinear relations, which have the lower-left matrix part interact with a suitably chosen block-diagonal part, are only semi-conservative in that  $G$-type and $J$-type elements are  preserved but appear in an operational context  showing  different direction of rotation. The second section begins with an examination of the immersion effect,  a tendency of LL$(G_{\mu,\nu}^{(p')})$ (for $G$-like lower-left part at order $p'\hspace*{-0.8mm}:=2^{n+1}\hspace*{-0.5mm}-\hspace*{-0.5mm}1$) to adopt  pecularities of LL$(J_{\mu,\nu}^{(p)})$: a semiprime  in LL$(J_{\mu,\nu}^{(15)})$ leads to a semiprime  in LL$(G_{\mu,\nu}^{(31)})$; diagonally repeated semiprimes  in LL$(J_{\mu,\nu}^{(31)})$ prompt diagonally repeated semiprimes in LL$(G_{\mu,\nu}^{(63)})$. After translating  opposing directions of rotation to differently directed paths that emerge when the ordinal number of the second factor of a semiprime with first factor 3, 7 or 11 is varied,  a rulebook is formulated which reveals that paths associated with LL$(J_{\mu,\nu}^{(p)})$ are not only aimed differently; they occur  five times as often as paths associated with LL$(G_{\mu,\nu}^{(p)})$. The coincidence of the frequency ratio ${\cal G}^{(p)}/{\cal J}^{(p)}$ with  the mass ratio of conventional matter over dark matter is commented upon at the close.
\end{abstract}
\maketitle

\section{Introduction -- Rational parafermi algebra}

\setlength{\fboxsep}{1mm}

 There has been  renewed interest in paraparticles as of lately \cite{Wang}. An attempt was also made to banish irrational numbers from quantum mechanics \cite{Palmer}. The ban would also affect types of parafermi algebra that are not based  entirely on rational numbers.  The original parafermi algebra, following  \cite{Green},  was directed at a
generalization of a spin-$(1/2)$ particle (fermion) to a spin-$(p/2)$ particle and contained irrational numbers:
\begin{equation}
b^{p+1}=\left(b^{+}\right)^{p+1}=0,\end{equation}
with the  $(p+1)\!\times\!(p+1)$
matrix representation for $b$, \begin{equation}
b_{\alpha\beta}=B_{\beta}\delta_{\alpha,\beta+1},\;\;(b^{+})_{\alpha\beta}=B_{\alpha}\delta_{\alpha+1,\beta},\qquad B_{\beta}=\sqrt{\beta(p-\beta+1)},\end{equation}
realizing the spin-$(p/2)$ representation\begin{equation}
\frac{1}{2}[b^{+},b]=\mathrm{diag}\textrm{\Large(}\frac{p}{2},\frac{p}{2}-1,\cdots,-\frac{p}{2}+1,-\frac{p}{2}\textrm{\Large)}\label{eq:b-spin}\end{equation}
 and the characteristic trilinear relations of parafermi algebra,\begin{equation}
[[b^{+},b],b]=-2b,\;\:[[b^{+},b],b^{+}]=2b^{+}.\label{eq:trilinear}\end{equation}
As was demonstrated in \cite{Merkel},  if the scope of  operators is restricted to Mersennian order $p=2^{n}-1\hspace*{2mm}n=1,2,\dots$, one can establish a Continued Kronecker Square-root process (CKS) 
\begin{equation} \sqrt{f\otimes I_2}, \hspace*{2mm}\sqrt{\sqrt{f\otimes I_2}\otimes I_2}, \hspace*{2mm}\sqrt{\sqrt{\sqrt{f\otimes I_2}\otimes I_2}\otimes I_2}, \hspace*{2mm}\dots
\end{equation}
starting with   $f\equiv f^{(1)}=\left(\begin{array}{cc}
0 & 0\\
1 & 0
\end{array}\right)$, that dispenses with irrational numbers altogether.  The question arises: Why continually switch to block matrices instead of using whatever is there? The CKS $-$which does exactly that$-$ does not deny that extracting the square root of  $f$ is impossible. The  Kronecker method of extracting the square root just offers a  way around this difficulty. The sequence begins with\,\footnote{\,omitting the minus sign is a venial sin  since we deal with a continued process}\begin{equation} \sqrt{f\otimes I_2}=
\left(\begin{array}{cccc}
0 & 0& 0& 0\\
0 & 0& 0& 0\\
1 & 0& 0& 0\\
0& 1& 0& 0\\
\end{array}\right)^{\hspace*{-1mm}1/2}=\left(\begin{array}{cccc}
0 & 0& 0& 0\\
1 & 0& 0& 0\\
0 & 1& 0& 0\\
y & 0& 1& 0\\
\end{array}\right)
\end{equation}
$-$where $y$ can be any number. In order to structure the CKS, we lay down that the choice $y=-1$ supplements the lower left block $c_3=\left(\begin{array}{cc}
0 & 1\\
-1 & 0
\end{array}\right)$ and initiates a CKS as adressed above, while
 the choice $y=1$  leads to a lower left block $c_2=\left(\begin{array}{cc}
0 & 1\\
1 & 0
\end{array}\right)$,   iniating another  CKS
\begin{equation} \sqrt{h\otimes I_2}, \hspace*{2mm}\sqrt{\sqrt{h\otimes I_2}\otimes I_2}, \hspace*{2mm}\sqrt{\sqrt{\sqrt{h\otimes I_2}\otimes I_2}\otimes I_2}, \hspace*{2mm}\dots
\end{equation} where $h=h^{(1)}= f^{(1)}=(c_2-c_3)/2$.%
\footnote{\ \  ${f}^{(1)}={h}^{(1)}=\left(\begin{array}{cc}
0 & 0\\
1 & 0
\end{array}\right)$,  ${I_2}\equiv\left(\begin{array}{cc}
1 & 0\\
0 & 1
\end{array}\right)\;\textrm{where} \;  c_{3}=\left(\begin{array}{cc}
0 & 1\\
-1 & 0
\end{array}\right)$, $c_{2}=\left(\begin{array}{cc}
0 & 1\\
1 & 0
\end{array}\right)\!\!\!$ \\

constitute the Clifford algebra $Cl(3)$
}%
\hspace*{1mm} Referring the reader to \cite{Merkel} for details, we  emphasize  that the relevant constructions (lower triangular  matrices with  row indices $\mu$ and column indices $\nu$)
\begin{eqnarray}\nonumber
 f^{(p)}=&I_{(p+1)/2}\;f^{(1)}+(G_{\mu\nu}^{(p)})\otimes c_3\\
 & & (\mu>\nu),\label{eq:rationalnews}\\ \nonumber
h^{(p)}=&I_{(p+1)/2}\;h^{(1)}+(J_{\mu\nu}^{(p)})\otimes c_2 \nonumber
\end{eqnarray}

\noindent  contain only natural numbers $G_{\mu\nu}^{(p)}$, $J_{\mu\nu}^{(p)}$ as  elements.

\noindent The two kinds of construction share a parafermi algebra that closely  resembles the Green's relations. The spin-$(p/2)$ representation is realized recursively:

\medskip{}

\noindent for order $p=1:\hspace*{3,5mm}s_1=\frac{1}{2}[c_2^T,c_2]=\mathrm{diag}(\frac{1}{2},\frac{-1}{2})$;

\medskip{}

\noindent
for order $p=3:\hspace*{3.5mm}s_3=I_2\otimes s_1+[c_2^T,c_2]\otimes I_2=\mathrm{diag}(\frac{3}{2},\frac{1}{2},\frac{-1}{2},\frac{-3}{2})$;

\medskip{}

\noindent
for order $p=7:\hspace*{3.5mm}s_7=I_2\otimes s_3+2\, [c_2^T,c_2]\otimes I_4=\mathrm{diag}(\frac{7}{2},\frac{5}{2},\,...\,,\frac{-5}{2},\frac{-7}{2})$;
\medskip{}

\noindent or generally,

\hspace{12mm} \begin{eqnarray} \nonumber
\hspace*{2mm}s_p=I_2\,\otimes\, s_{(p-1)/2}+{\textstyle\frac{p+1}{4}} \left[c_2^T,c_2\right]\,\otimes \,I_{(p+1)/2}\\ (p=3,7,\dots)\label{eq:b-spin}\\ =\,\mathrm{diag}\textrm{\Large(}\frac{p}{2},\frac{p}{2}-1,\cdots,-\frac{p}{2}+1,-\frac{p}{2}\textrm{\Large).}\nonumber
\end{eqnarray}

\smallskip{}
\noindent In the spin representation,  recursion manifests as a uniform doubling of the diagonal's reach. 
But as far as  $f^{(p)}$ and  $h^{(p)}$ are concerned: due to their triangular shape, recursion implies that the greatest increase in information occurs on the lower-left quadrant (LL); existing information is (redundantly) distributed across the upper and lower triangles, plus a duplicate in LL's upper-right subquadrant.

\noindent The analogues to Smith's trilinear relations (Eq.\,\ref{eq:trilinear}) reveal that the information collections LL$(G_{\mu,\nu})$ and  LL$(J_{\mu,\nu})$, as fixated to spin $c_k$ along the $k$-th coordinate $(k=3,2,)$ have opposite direction of rotation, and, in that sense, are only semi-conservative:
where $f_{\textrm{LR}}^{(p)}
=h_{\textrm{LR}}^{(p)}=f\otimes I_{(p+1)/4} \;\;(p=3,7,\dots,)$

%\medskip{}
\begin{eqnarray}\nonumber
\left[\;\left[(f_{\textrm{LR}}^{(p)})^{T},\;\textrm{LL}(G_{\mu\nu}^{(p)})\otimes c_3\right],\,f_{\textrm{LR}}^{(p)}\;\right]\;=\;2\,\textrm{LL}(G_{\mu\nu}^{(p)})\otimes f,
\nonumber\\
\left[\;\left[(f_{\textrm{LR}}^{(p)})^{T},\;\textrm{LL}(G_{\mu\nu}^{(p)})\otimes c_3\right],\,(f_{\textrm{LR}}^{(p)})^{T}\right]\;=-2\,\textrm{LL}(G_{\mu\nu}^{(p)})\otimes f^T,\nonumber\\
\\\label{eq:trilinearrational}
\left[\;\left[(h_{\textrm{LR}}^{(p)})^{T},\;\textrm{LL}(J_{\mu\nu}^{(p)})\otimes c_2\right],\,h_{\textrm{LR}}^{(p)}\;\right]\;=-2\,\textrm{LL}(J_{\mu\nu}^{(p)})\otimes f,
\nonumber\\
\left[\;\left[(h_{\textrm{LR}}^{(p)})^{T},\;\textrm{LL}(J_{\mu\nu}^{(p)})\otimes c_2\right],\,(h_{\textrm{LR}}^{(p)})^{T}\right]\;=\;2\,\textrm{LL}(J_{\mu\nu}^{(p)})\otimes f^T.\nonumber
\end{eqnarray}

\noindent  

\noindent For order $p=31$, the lower-left quadrants (with their upper-right subquadrants framed) look like

\begin{figure}[H]
\caption{\label{fig:catal-rep} Order-31 arrangement of  information on 
${f}^{(31)}$ and ${h}^{(31)}$}
\[
\mathrm{LL}\,(G_{\mu\nu}^{(31)})=\qquad\qquad\qquad\qquad
\]
$\quad\left(\begin{array}{c}
\\
\\
\\
\\
\\
\\
\\
\\
\end{array}\right.$\hspace{-0.3cm}%
\begin{tabular}{cccc|cccc}
$\underline{429}$ & 155 & 43 & 19 & $\underline{5}$ & \multicolumn{1}{c|}{3} & $\underline{1}$ &|$\underline{\,1\;\,}$\tabularnewline
1275 & $\underline{429}$ & 115 & 43 & 11 & \multicolumn{1}{c|}{$\underline{5}$} & 1 & $\underline{1}$\tabularnewline
\cline{7-8} 
4819 & 1595 & $\underline{429}$ & 155 & 41 & 17 & $\underline{5}$ & 3\tabularnewline
15067 & 4819 & 1275 & $\underline{429}$ & 113 & 41 & 11 & $\underline{5}$\tabularnewline
\cline{5-8} 
58781 & 18627 & 4905 & \multicolumn{1}{c}{1633} & $\underline{429}$ & 155 & 43 & 19\tabularnewline
189371 & 58781 & 15297 & \multicolumn{1}{c}{4905} & 1275 & $\underline{429}$ & 115 & 43\tabularnewline
737953 & 227089 & 58781 & \multicolumn{1}{c}{18627} & 4819 & 1595 & $\underline{429}$ & 155\tabularnewline
2430289 & 737953 & 189371 & \multicolumn{1}{c}{58781} & 15067 & 4819 & 1275 & $\underline{429}$\tabularnewline
\end{tabular}\hspace{-0.2cm}$\left.\begin{array}{c}
\\
\\
\\
\\
\\
\\
\\
\\
\end{array}\right)$\medskip{}
\[
\mathrm{LL}\,(J_{\mu\nu}^{(31)})=\qquad\qquad\qquad\qquad
\]
 $\left(\begin{array}{c}
\\
\\
\\
\\
\\
\\
\\
\\
\end{array}\right.$\hspace{-0.3cm}%
\begin{tabular}{cccc|cccc}
$-\underline{429}$ & 117 & -41 & 13 & $-$5\hspace*{-2mm}\color{red}{\ul{\,\,\,\,}}& \multicolumn{1}{c|}{1} & $-\underline{1}$ & |$\underline{\;1\;\,}$\tabularnewline
1547 & $-\underline{429}$ & 143 & -41 & 15 & \multicolumn{1}{c|}{$-5$\hspace*{-2mm}\color{red}{\ul{\,\,\,\,}}} & 3 & $-\underline{1}\,$\tabularnewline
\cline{7-8} 
-4903 & 1343 & $-\underline{429}$ & 117 & -43 & 15 & $-$5\hspace*{-2mm}\color{red}{\ul{\,\,\,\,}}
& 1\tabularnewline
18269 & -4903 & 1547 & $-429$ & 149 & -43 & 15 & $-$5\hspace*{-2mm}\color{red}{\ul{\,\,\,\,}}\tabularnewline
\cline{5-8} 
-58791\hspace*{-9mm}\color{red}{\ul{\quad\quad\,\,\,}} & 15547 & -4823 & \multicolumn{1}{c}{1319} & $-429$ & 117 & -41 & 13\tabularnewline
223573 & -58791\hspace*{-9mm}\color{red}{\ul{\quad\quad\,\,\,}} & 17989 & \multicolumn{1}{c}{-4823} & 1547 & $-429$ & 143 & -41\tabularnewline
-747765 & 194993 & -58791\hspace*{-9mm}\color{red}{\ul{\quad\quad\,\,\,}} & \multicolumn{1}{c}{15547} & -4903 & 1343 & $-\underline{429}$ & 117\tabularnewline
2886235 & -747765 & 223573 & \multicolumn{1}{c}{-58791\hspace*{-9mm}\color{red}{\ul{\quad\quad\,\,\,}}} & 18269 & -4903 & 1547 & $-\underline{429}$\tabularnewline
\end{tabular}\hspace{-0.2cm}$\left.\begin{array}{c}
\\
\\
\\
\\
\\
\\
\\
\\
\end{array}\right)$
\end{figure}

\smallskip{}

\noindent One of the unsolved questions is  why the emergent numbers appear to be  governed by Catalan numbers$-$when we zoom in on the upper-right subquadrant of $\mathrm{LL}\,(G_{\mu\nu}^{(31)})$, for example, we observe how they turn up alongside the secondary diagonal in a non-interfering way:
\begin{eqnarray} \nonumber
\underline{5}=C_3,\\
3+11=C_4, \\ \nonumber
1+41=C_5, \\ \nonumber
 1+1+17+113=C_6,
\end{eqnarray}
and the next round starts with $\underline{429}=C_7$.
 Catalan numbers grow approximately like 
\begin{equation} C_k\sim\,\frac{2}{\sqrt \pi}\,\frac{(\sqrt2)^\ell}{\ell^{\,3/2}}\,,\end{equation}
where $k=\ell/4-1$. 

\smallskip{}
\noindent Two diagonals in the subquadrants of LL$(J_{\mu\nu}^{(31)})$ have been underlined in red. In the upper-right, there reside copies of $-5$. In the lower-left, copies of $-58791$; when taking $k$ modulo 8 and the  elements by their absolute value,  $\vert-58791\vert$ turns  the Catalan resonance $C_3+C_{11}$, while the corresponding element for  LL$(G_{\mu\nu}^{(31)})$, 58781, the `Catalan beat' $C_{11}-C_3$\,.

%Both relationships possess an inherent recursivity. The spin-$p/2$ representation is overtly recursive relative to the preceding epoch; 
%s is the trilinear relations' recursion which resides in the arrangement of pasts within $(G_{\mu\nu}^{(p)}\otimes c_3)$; it shows up by its staggeredness:  the upper left (or lower right) arrangement pserves the stages $n-1, n-2, ... , 1$ ($n=\log_2(p+1)$), while the lower left arrangement refers to the present, $n$, followed  by the stages $n-1, n-2, ... , 2$ (framed subquadrants). .
\smallskip{}
\noindent To remain at  order 31, the new information brought to light by  CKS  can be compressed into the 18-tuples
\begin{eqnarray*}
(G_{\rho}^{(31)}) & = & (19,43,115,155,\underline{429},1275,1595,1633,4819,4905,\\
 &  & \:15067,15297,18627,58781,189371,227089,737953,2430289)
\end{eqnarray*}
and
\begin{eqnarray*}
(J_{\rho}^{(31)} ) & = & (13,-41,117,143,-\underline{429},1319,1343,1547,-4823,-4903,\\
 &  & 15547,17989,18269,-58791\textrm{\hspace*{-8.5mm}\textcolor{red}{\ul{\qquad\,\,}}},194993,223573,-747765,2886235).
\end{eqnarray*}

%\noindent Not all
%of the quadrants
%need to be considered because the subquadrants UL(LL($\cdot$)) =
%LL(UL($\cdot$)) = LL(LR($\cdot$)) just reproduce order-$(p-1)/2$ 
%quantities. New quantities are named base numbers. As shown in Fig. \ref{fig:catal-rep}, order-31 base numbers
 %can be extracted from the non-UR(LL($\cdot$)) parts
%of quadrants $\textrm{LL}(G_{\mu\nu}^{(31)})$ and $\textrm{LL}(J_{\mu\nu}^{(31)})$:
%
%\noindent 

%\newpage
%\noindent from which the outer nodes of 18-cube complexes with boundary
%labels $\Gamma_{x}$ and $\chi_{x}$ are respectively formed. 

%\vspace*{2mm}
%\newpage
\noindent For later use, let's list  past news in the order it appears and weathers away:

\bigskip{}
$\begin{array}{ll}(G_{\rho}^{(3)})=\left(1\right), &(J_{\rho}^{(3)})=\left(1\right),\\
                           (G_{\rho}^{(7)})=\left(1\right), &(J_{\rho}^{(7)})=\left(-1,3\right),\\       
                           (G_{\rho}^{(15)})=\left(3,5,11,17,41,113\right), &(J_{\rho}^{(15)})=\left(1,-5,15,-43,149\right);\\                         
\end{array}$

\bigskip{}
\noindent the order-31 tuples have been shown above; for order 63,
at least one tuple, with  54 members, is known,
\[
\begin{array}{c}
(G_{\rho}^{(63)})=\\
\\
\left(\begin{array}[b]{r}
\\
\\
\end{array}\right.
\end{array}\!\!\!\!\!\!\!\!\!\!\!\!\!\!\!\!\!\!\!\!\!\!\!\!\!\!\!\!\!\!\begin{array}[t]{rrr}
1,\!635\,, & 4,\!907\,, & 15,\!299\,,\\
18,\!633\,, & 58,\!791\textrm{\hspace*{-9.5mm}\textcolor{red}{\ul{\qquad\,\,\,}}}\,, & 189,\!393\,,\\
227,\!123\,, & 230,\!395\,, & 738,\!035\,,\\
747,\!851\,, & 2,\!430,\!515\,, & 2,\!461,\!115\,,\\
2,\!926,\!323\,, & \underline{9,\!694,\!845}\,, & 32,\!431,\!347\,,\\
38,\!331,\!419\,, & 38,\!792,\!251\,, & 38,\!795,\!521\,,\\
128,\!896,\!939\,, & 130,\!392,\!731\,, & 130,\!402,\!545\,,\\
436,\!615,\!771\,, & 441,\!538,\!235\,, & 441,\!568,\!833\,,\\
519,\!693,\!033\,, & 519,7\!30,\!299\,, & 1,\!767,\!204,\!399\,,\\
1,\!767,\!321,\!981\,, & 6,\!044,\!219,\!361\,, & 6,\!044,\!598,\!147\,,\\
7,\!090,\!735,\!179\,, & 7,\!091,\!189,\!425\,, & 7,\!168,\!783,\!827\,,\\
24,\!335,\!131,\!347\,, & 24,\!336,\!607,\!417\,, & 24,\!597,\!422,\!507\,,\\
83,\!908,\!823,\!403\,, & 83,\!913,\!684,\!433\,, & 84,\!796,\!853,\!171\,,\\
99,\!228,\!108,\!067\,, & 343,\!059,\!613,\!221\,, & 1,\!190,\!676,\!037,\!827\,,\\
1,\!390,\!379,\!604,\!203\,, & 1,\!404,\!717,\!639,\!489\,, & 4,\!837,\!348,\!974,\!083\,,\\
4,\!886,\!545,\!335,\!065\,, & 16,\!885,\!007,\!814,\!155\,, & 17,\!054,\!606,\!505,\!569\,,\\
19,\!881,\!172,\!597,\!035\,, & 69,\!531,\!783,\!535,\!237\,, & 243,\!860,\!214,\!616,\!867\,,\\
283,\!858,\!869,\!110,\!417\,, & 997,\!331,\!203,\!563,\!441\,, & 3,\!512,\!576,\!820,\!924,\!177
\end{array}\begin{array}[t]{r}
\\
\\
\\
\\
\\
\\
\\
\\
\\
\\
\\
\\
\\
\\
\\
\left.\!\!\!\!\!\!\begin{array}{r}
\\
\\
\\
\end{array}\right),
\end{array}
\]
\noindent where, due to space constraints,  only an excerpt  of the   matrix (5 rows, each truncated)  it extracts new information from can be shown:

\smallskip{}

$\mathrm{LL}\,(G_{\mu\nu}^{(63)})=\begin{array}{cc}
\\
\\
\\
\end{array}$ 

\noindent

\hspace*{-0.65cm}\begin{minipage}[t]{1\columnwidth}%
$\left(\begin{array}{rrrrrrrrr}
\scriptstyle{\underline{9694845}} & \scriptstyle{2926323} & \scriptstyle{747851} & \scriptstyle{230395} & \scriptstyle{58791}\hspace*{-7mm}\color{red}{\underline{\qquad}} & \scriptstyle{18633} & \scriptstyle{4907} & \scriptstyle{1635} & \ldots\\
\scriptstyle{32431347} & \scriptstyle{\underline{9694845}} & \scriptstyle{2461115} & \scriptstyle{747851} & \scriptstyle{189393} & \scriptstyle{58791}\hspace*{-7mm}\color{red}{\underline{\qquad}} & \scriptstyle{15299} & \scriptstyle{4907} & \ldots\\
\scriptstyle{128896939} & \scriptstyle{38331} & \scriptstyle{\underline{9694845}} & \scriptstyle{2926323} & \scriptstyle{738035} & \scriptstyle{227123} & \scriptstyle{58791}\hspace*{-7mm}\color{red}{\underline{\qquad}} & \scriptstyle{18633} & \ldots\\
\scriptstyle{436615771} & \scriptstyle{128896939} & \scriptstyle{32431347} & \scriptstyle{\underline{9694845}} & \scriptstyle{2430515} & \scriptstyle{738035} & \scriptstyle{189393} & \scriptstyle{58791}\hspace*{-7mm}\color{red}{\underline{\qquad}} & \ldots\\
\scriptstyle{1767204399} & \scriptstyle{519693033} & \scriptstyle{130392731} & \scriptstyle{38792251} & \scriptstyle{\underline{9694845}} & \scriptstyle{2926323} & \scriptstyle{747851} & \scriptstyle{230395} & \ldots\\
 &  &  & &  &  \ldots &  &  &   
\end{array}\right)$%
\end{minipage}

%\noindent 

%\[
%\]
{\small$(\mu=17,\dots,32;\,\nu=1,\dots,16.)$}

\medskip{}
% \newpage
\noindent The  tuple members can be grouped as follows:

\[
 \begin{array}{l}
\!\!\!\!\!\!\!\!\!\!\!\!\!\!\!\!\!\!(3^{\log_{2}(p+1)-3}+1)/2 \quad\textrm{tuple members}\;
\leq C_{q} \\
\\
\!\!\!\!\!\!\!\!\!\!\!\!\!\!\!\!\!\!(3^{\log_{2}(p+1)-2}-1)/2 \quad\textrm{tuple members}\; 
>  C_{q} \\
\\
\!\!\!\!\!\!\!\!\!\!\!(q=(p-3)/4,\quad  p>15). \\
\end{array}
\]

\smallskip{}
\noindent It is striking that little to no newness appears up until order 15. One possible explanation for this behavior could be linked to another form of `semi-'ness, which we will  discuss now.
\smallskip{}
%\noindent Only if we were to encode the terms of a Charles Dodgson universe in their entirety would these tuple members apply.)
%\newpage
\section{The role of $G_\rho^{(p)}$ and $J_\rho^{(p)}$, semiprime with first factor 3, 7 or 11}

%\medskip{}
%One of the basic identities for semiprimes is given by
%\begin{equation}\pi(m)\times\pi(m+1)\equiv(m+1)\times\pi(m+1)\;\;\big(\textrm{mod}\;\pi(m+1)\big)\label{eq:trivial}
%\end{equation} where $\pi(m)$ denotes the prime number with ordinal number $m$. It is trivially satisfied by $3\times5$. More remarkable is the fact that it equals the order it occurs with: $J_3^{(15)}=15$. That order is also bound up with another trivial instance of Eq.\,\ref{eq:trivial}\,: $\#\,G_\rho^{(15)}=2\times3$.
%Despite their apparently slow development into  newness, the constructions have surprises in store. In the  picture they provide, there's more than one ``bang,'' so to say. At order 7, we observe  that one element from LL$(J_{\mu\nu}^{(7)})$  breaks loose from identity: $J_{2,1}^{(7)}=3$.  A second ``bang'' occurs at order 15, which then  encompasses elements from both LL$(J_{\mu\nu}^{(15)})$ and LL$(G_{\mu\nu}^{(15)})$. It's to be expected that processes like these leave their imprint on their specific surroundings. During the  parafermionic evolution, there  exist one or more information  paths toward orders 7,\,15. We found   that these paths take the form of a sequence of distinct pairs of semiprimes. The way they choose follows a rulebook that will be presented next.  
In \cite{Merkel}, we discussed how $G_{\mu\nu}^{(p)}$,$\,J_{\mu\nu}^{(p)}$ influence their Mersennian successors,  $G_{\mu\nu}^{(p')}$,$\,J_{\mu\nu}^{(p')}$, an influence which we called  interordinality. A special side of the phenomenon is the immersion effect$-$when new details emerge in $\,J_{\mu\nu}^{(p)}$, then $G_{\mu\nu}^{(p')}$ immerses into the new `milieu' and adopts those details as if they were a part of itself.  At order 7,  one element from LL$(J_{\mu\nu}^{(7)})$ comes up with such a detail: it breaks loose from identity, $J_{2,1}^{(7)}=3$.  Immediately after, six elements  from  LL$(G_{\mu\nu}^{(15)})$ break loose from identity. At that order, LL$(J_{\mu\nu}^{(15)})$ has something novel  in store, a semiprime: $J_{3}^{(15)}=3\times5$; it takes just one step further,  and  LL$(G_{\mu\nu}^{(31)})$ has   a semiprime of its own:
$G_{10}^{(31]}=3\times5099$. No surprise, then, that an entire diagonal of LL(LL$(J_{\mu\nu}^{(31)})$) finds its counterpart in UR(UL(LL$(G_{\mu\nu}^{(63)})))$, containing, except a minus sign, the same  copies of 58791. Returning to $J_{3}^{(15)}=3\times5$, we note that
it is not just any semiprime$-$it has a value that coincides with the order it occurs on. Moreover,  it is a prime-indexed prime (PIP) with iteration index 2: \begin{equation}3\times5=3\pi(3)=3\pi(\pi(2)),\end{equation}which  means it might hold a special place in some larger scheme. There exists one identity for consecutive primes like 3 and 5,
\begin{equation}\pi(m)\times\pi(m+1)\equiv(m+1)\times\pi(m+1)\;\;\big(\textrm{mod}\;\pi(m+1)\big)\label{eq:trivial}
\end{equation} where $\pi(m)$ denotes the prime number\,\footnote{\,not the circle constant, which can't interfere in the present context} with ordinal number $m$.
%It seems logical, then, to start with a pair $(\pi(m),\pi(m+1))$. In fact, the first departure from unity occurs at order 7, and the LL block the deviating element  springs from has determinant 2, so $m+\log_2\,2$  would justify the  pair's second-member form.
%Now 7 and 15
% are both Mersennian orders $(p=2^n-1)$ which can be considered a sign-reversed involution: If $p(x)$ is that function, then $-p(-p(2^n)) =2^n$; this suggests using a pseudo-involutive form for a pair's second member, \begin{equation}\big(\pi(m)\times\pi(m+1),\pi(\pi(m)\times\pi(m+1)\big).\label{eq:plus1}\end{equation} Applied to  $(J_{\nu+1,\nu}^{(15)})=3\times5$, this gives  \begin{equation}(3\times5,\pi(3\times7)),\end{equation} with \begin{equation}(3\times7,\pi(3\times11))\label{eq:continuation}\end{equation} as its successor. 
 %The second bang offers one further increment: The LL$(J_{\mu\nu}^{(15)})$  superblock has four blocks, with determinants $\det \textrm{UL}=\det \textrm{LR}=10,\det \textrm{UR}=2$ and $\det \textrm{LL}=386$. Their sum modulo 16 gives 8, furnishing a new increment $\log_2 8=3$. so we have to look for a new pairing \begin{equation}\big(\pi_a\times\pi(m),\pi(\pi_b\times\pi(m+3)\big).\end{equation} 
% Since this is the second bang's environ, a  match has  to be accomplished with the former pair's continuation, \begin{equation}(7\times3,\pi(3\times11)),\label{eq:contalter}\end{equation}
%giving $\pi_a=7,\pi_b=3$; the new pair, however,  leads  to a path that differs from that produced by pair (\ref{eq:continuation}).
The second factor of $G_{10}^{(31]}$, though, is a semiprime-indexed prime (SPIP) with left-iteration index 1 and  right-iteration index 3.  \begin{equation}5099=\pi(3\times227)=\pi(3\pi(7\times7))=\pi(3\pi(7\pi(2\times2))).\label{eq:iteration}\end{equation}
It turns out there  exists another  identity for this type of semiprime: Where $\pi(5)=11,
\pi(10)=29,\pi(49)=227$, its first three instantiations are \begin{equation}\begin{array}{cccc}\pi(3\times11)&\equiv&11&\hspace*{-2mm}(\textrm{mod}\;3\times7),\\
\pi(3\times29)&\equiv&29&\hspace*{-2mm}(\textrm{mod}\;3\times7),\\
\pi(3\times227)&\equiv&227&\hspace*{-2mm}(\textrm{mod}\;3\times7).\label{eq:involfurther}\end{array}\end{equation}
So the immersion effect is best characterized by a pair \begin{equation}\big(3\times\pi(m),\pi(3\times\pi(m+1)\big)\label{eq:basicpair},\end{equation} with first instance \begin{equation}(3\times5,\pi(3\times7))\end{equation} and first continuation \begin{equation}(3\times7,\pi(3\times11))\label{eq:nuclear}.\end{equation} It's also clear that one can swap the first pair's factors:\begin{equation}(7\times3,\pi(3\times11))\label{eq:swap37}\end{equation}  provided we adapt  ordinal numbers accordingly:  \begin{equation}\big(7\times\pi(m),\pi(3\times\pi(m+3)\big),\end{equation} which, of course, leads to a different sequence of pairs as $m$ gets varied. When we swap the factors of the second pair, too: \begin{equation}(7\times3,\pi(11\times3)),\end{equation} according to  \begin{equation}\big(7\times\pi(m),\pi(11\times\pi(m)\big),\end{equation} the sequence of pairs will vary once again.

\smallskip{}
\noindent Those variations should suffice$-$as a minimal rulebook  for testing the subset of semiprimes with first factors 3, 7 and 11 as they turn up for order 31:

\smallskip{
}
\begin{center}$\begin{array}{rcl}G_{10}^{(31]}=&15297=&\hspace*{1.8mm}3\times5099,\\ 
                                      J_4^{(31)}=&\;143\;\;\,=&11\times13, \\
                                      J_9^{(31)}=&15547=&\hspace*{1.8mm}7\times2221%(=7\pi(331))=7\pi(\pi(67))=7\pi(\pi(\pi(19)))=7\pi(\pi(\pi(\pi(8))))
                                      ,\\
                                      \vert J_{12}^{(31)}\vert=&58791=&\hspace*{1.8mm}3\times19597 %(=3\pi(2223)).
                                      \end{array}$
                                      \end{center}

%\noindent Somewhat surprisingy, there's  only one $G$-type first-factor occurrence,  $G_{\rho}^{(31]}$, which is bound up with a  pseudo-involutive game taken even further via identity (\ref{eq:involfurther}): \begin{equation}G_{10}^{(31]}=15297=3\times5099=3\,\pi(681)=3\,\pi(3\times227)=3\times\pi(3\,\pi(49)).\end{equation}
%\noindent$3\times5$.Rule I,Table \ref{tab:anterograde-pathway} however encompasses $3\times7$, too
%\setul{1pt}{1.0pt}
\begin{table}[H]
\caption{\label{tab:anterograde-pathway} An anterograde path  which leads from  $3\times2$  to $3\times5$ and, from that point on,  to $3\times7$ where Rule I says that  $3\pi(m)$  ($m=1,2,3,4$) is to be paired with 
 $\pi(3\pi(m+1))$;}
\vspace{0.3cm}
\qquad{}\qquad{}\qquad{}
\begin{tabular}{|c|c|c|}
\hline 
 & \multicolumn{2}{c|}{$\begin{array}{ccccc}
\\
\\
\\
\end{array}$Rule I}\tabularnewline
ordinal number $m$& $3\pi(m)$ & $\pi(3\pi(m+1))$\tabularnewline
\hline 
\hline 
1 & $\textcolor{black}{3\times2}$& \hspace*{-8.mm}$\textcolor{lightgray}{\pi(3\times3)}$
\tabularnewline
\hline 
2\hspace*{0mm}& $\textcolor{black}{3\times3}$&  \hspace*{-8.mm}$\textcolor{lightgray}{\pi(3\times5)}$\tabularnewline
\hline 
\hspace*{12.5mm}\textcolor{black}{3}\hspace*{5mm}$J_3^{(15)}$ &\textcolor{blue}{$3\times5$}\hspace*{-8mm}\color{red}{\ul{\quad\quad\,\,}} 
 & \hspace*{-7.5mm}$\textcolor{lightgray}{\pi(3\times7)}$\tabularnewline
\hline 
4 &$\textcolor{red}{3\times7}$  & \hspace*{-5,5mm}$\textcolor{lightgray}{\pi(3\times11)}$  \tabularnewline
\hline 
%\ul{soeinquatsch}
\end{tabular}
\end{table}

\begin{table}[H]
\caption{\label{tab:retrograde-pathway3} Two  paths, one anterograde  for $58791=\vert J_{12}^{(31)}\vert=3\times19597$, and another retrograde for $G_{10}^{(31)}=3\times5099$  ending at $\pi(3\times11)$, where  Rule I says that  $3\pi(m)$  ($m\geq4$) is to be paired with 
 $\pi(3\pi(m+1))$.}
\vspace{0.1cm}
\qquad{}\qquad{}\qquad{}%
\begin{tabular}{|c|c|c|}
\hline 
 & \multicolumn{2}{c|}{$\begin{array}{ccccc}
\\
\\
\\
\end{array}$Rule I}\tabularnewline
ordinal number $m$& $3\pi(m)$ & $\pi(3\pi(m+1))$\tabularnewline
\hline 
\hline 
\hspace*{-1.45cm}\vdots&\vdots&\vdots\tabularnewline
\hline 
$\textcolor{black}{2224\hspace*{1.45cm}}$& $\textcolor{black}{3\times19603}$&$\textcolor{lightgray}{\pi(3\times19609)}$\tabularnewline
\hline 
2223\hspace*{5mm}$\vert J_{12}^{(31)}\vert$ & $\textcolor{blue}{3\times19597}$\hspace*{-15mm}\color{red}{\ul{\qquad\qquad\,\,}} & $\textcolor{lightgray}{\pi(3\times19603)}$\tabularnewline
\hline 

\hspace*{-1.45cm}\vdots&\vdots&\vdots\tabularnewline
\hline 
\hspace*{-1.45cm}680& $\textcolor{lightgray}{3\times5087}$& \hspace*{12mm}$\pi(\textcolor{blue}{3\times5099})$\hspace*{5mm}$G_{10}^{(31)}$ \tabularnewline
\hline 
\hspace*{-1.45cm}\vdots&\vdots&\vdots\tabularnewline
\hline 
%\textcolor{black}{226\hspace*{1.45cm}}$ &$\textcolor{lightgray}{3\times1429}$  &$\textcolor{blue}{\pi(3\times681)}$\tabularnewline
%\hline 
 5\hspace*{1.52cm} &$\textcolor{lightgray}{3\times11}$  &$\textcolor{black}{\pi(3\times13)}$  \tabularnewline
\hline 
4\hspace*{1.52cm} &$\textcolor{lightgray}{3\times7}$  &$\pi(\textcolor{red}{3\times11})$  \tabularnewline
\hline 

\end{tabular}
\end{table}

\begin{table}[H]
\caption{\label{tab:retrograde-pathway7} Two anterograde paths for $\pi(\vert J_{12}^{(31)}\vert)=\pi(58791)=\pi(3\times19597)$, $J_{9}^{(31)}=7\times2221$, % leading to $\pi(3\times11)=137$ and  $7\times3$,
 where Rule III says that  $7\pi(m)$  ($m\ge2$) is to be paired with 
 $\pi(3\pi(m+3))$;}
\vspace{0.5cm}
\qquad{}\quad{}%
\begin{tabular}{|c|c|c|}
\hline 
 & \multicolumn{2}{c|}{$\begin{array}{ccccc}
\\
\\
\\
\end{array}$Rule III}\tabularnewline
ordinal number $m$& $7\pi(m)$ & $\pi(3\pi(m+3))$\tabularnewline
\hline 
\hline 
\hspace*{-1.45cm}\vdots&\vdots&\vdots\tabularnewline
\hline 
\textcolor{black}{2221}\hspace*{14mm} &$\textcolor{black}{7\times19577}$  &$\textcolor{black}{\pi(3\times19603)}$\tabularnewline
\hline 
2220 \hspace*{13mm}& $\textcolor{black}{7\times19571}$& \hspace*{14mm}$\pi(\textcolor{blue}{3\times19597})$\hspace*{-16mm}\color{red}{\ul{\qquad\qquad\,\,}}\hspace*{5mm}\textcolor{black}{$\vert J_{12}^{(31)}\vert$}\tabularnewline
\hline 
%$\textcolor{black}{2219\hspace*{1.45cm}}$ &$\textcolor{lightgray}{7\times19559}$  &$\textcolor{black}{\pi(3\times19583)}$\tabularnewline
%\hline 
\hspace*{-1.45cm}\vdots&\vdots&\vdots\tabularnewline
\hline 
331\hspace*{7mm}$J_{9}^{(31)}$ & $\textcolor{blue}{7\times2221}$& $\textcolor{lightgray}{\pi(3\times2243)}$\tabularnewline
\hline 

3\hspace*{1.52cm} &$\textcolor{lightgray}{7\times5}$  &$\textcolor{lightgray}{\pi(3\times13)}$  \tabularnewline
\hline 
2\hspace*{1.52cm} &$\textcolor{lightgray}{7\times3}$  &$\textcolor{lightgray}{\pi(3\times11)}$  \tabularnewline
\hline 
\end{tabular}
\end{table}
%$\mathcal G/\mathcal J=(0+1+0)/(1+2+2)=1/5$;

%\noindent and for order 63, $\mathcal G/\mathcal J=(4+4+0)/(?+?+?)=?$ (see Tables 2-4.)

\begin{table}[H]
\caption{\label{tab:retrograde-pathway11} Two anterograde paths for  $J_{9}^{(31)}=7\times2221$ and $\pi(J_{4}^{(31)})=\pi(143)=\pi(11\times13)$, % towards $\pi(11\times3)=137$
 where Rule II says that  $7\pi(m)$  ($m\ge2$) is to be paired with 
 $\pi(11\pi(m))$;}
\vspace{0.5cm}
\qquad{}\qquad{}\qquad{}%
\begin{tabular}{|c|c|c|}
\hline 
 & \multicolumn{2}{c|}{$\begin{array}{ccccc}
\\
\\
\\
\end{array}$Rule II}\tabularnewline
ordinal number $m$& $7\pi(m)$ & $\pi(11\pi(m))\hspace*{11mm}$\tabularnewline
\hline 
\hline 
$\vdots\hspace*{13mm}$ &$\vdots$  & $\vdots$\hspace*{12mm} \tabularnewline
\hline 
331\hspace*{5mm}$J_{9}^{(31)}$ & $7\times\textcolor{blue}{2221}$& $\textcolor{black}{\pi(11\times2221)}\hspace*{8mm}$\tabularnewline
\hline 
330\hspace*{13mm} & $\textcolor{lightgray}{7\times2213}$& $\textcolor{black}{\pi(11\times2213)}\hspace*{8mm}$\tabularnewline
\hline 
$\vdots\hspace*{13mm}$ &$\vdots$  & $\vdots$\hspace*{12mm} \tabularnewline
\hline 
$\textcolor{black}{6}\hspace*{13mm}$ & $\textcolor{lightgray}{7\times13}$ & $\pi(\textcolor{blue}{11\times13})$\hspace*{5mm}$J_{4}^{(31)}$ \tabularnewline
\hline 
%$\textcolor{black}{5}\hspace*{13mm}$ & $\textcolor{lightgray}{7\times11}$& $\textcolor{lightgray}{\pi(11\times11)}$\hspace*{12mm} \tabularnewline
%\hline 
$\vdots\hspace*{13mm}$ &$\vdots$  & $\vdots$\hspace*{12mm} \tabularnewline
\hline 
2\hspace*{13mm} & $\textcolor{lightgray}{7\times3}$ &$\textcolor{lightgray}{\pi(11\times3)}$\hspace*{13mm} \tabularnewline 
\hline 
\end{tabular}
\end{table}

\noindent %The notion  that Tables \ref{tab:anterograde-pathway}-\ref{tab:retrograde-pathway11} urge upon us  is that there's a commutation relation at play: A semiprime with first factor 3 can only run for a Rule III pathway test if there's a semiprime with first factor 7 opting likewise. Symbolically: 3,7 only if 7,3. If not, Rule I must suffice. In the same vein, 7,11 only if 11,7 for Rule II pathway tests. 
Let's restate in formal terms what Rules I\,-\,III and Tables \ref{tab:anterograde-pathway}\,-\,\ref{tab:retrograde-pathway11} urge upon us.

\noindent The basic operation consists of forming a pair \begin{eqnarray}\big(\textrm{F}_1\,\pi(m),\pi(\textrm{F}_2\,\pi(m+\kappa))\big)\hspace*{3cm}\\ \nonumber\qquad\;\;\;
(\textrm{F}_1,\textrm{F}_2)\in\{(3,3),(7,3),(7,11)\},\; m\in\mathbb{N},\;\; \kappa\in\{0,1,3\}.\end{eqnarray}

\noindent If feasible, one can combine this with the adjoint operation \begin{equation}\big(\textrm{F}_1\,\pi(m-\kappa),\pi(\textrm{F}_2\,\pi(m))\big),\quad\;\;\,
m\in\mathbb{N}\backslash\{1,..,\kappa+1\}.\end{equation}.

\noindent The total number of uses of   semiprimes $G_\rho^{(p)},J_\rho^{(p)}$ is characterized by the proportion ${\cal G}^{(p)}/{\cal J}^{(p)}$, and the individual uses under Rules I\,-\,III, by
\begin{equation}\textrm{\# of uses} \,:=\;\;2-\delta_{\textrm{F}_1,\textrm{F}_2}\,.\label{eq:uses}\end{equation}

\noindent  
% (3+12,3+10),(7+9+1,3+12+0),(7+9+1,11+4+1)
The distribution of $G_\rho^{(p)}$,\,$J_\rho^{(p)}$ is carried out   in a parity-conserving way, \begin{equation}\begin{array}{cccc}\rho_1+\rho_2+(2-\delta_{\textrm{F}_1,\textrm{F}_2})&\equiv& \textrm{F}_1\times\textrm{F}_2&\hspace*{-2mm}\big(\textrm{mod}\;2\big).
\end{array}\end{equation}
 \begin{equation}\begin{array}{cccc}\multicolumn{4}{l}{\hspace*{-3.2cm}\textrm {\small Rule I,\,\,\,\,  Table \ref{tab:retrograde-pathway3}\,:}}\\
\hspace*{-1mm}12+10+1&\hspace*{-1mm}\equiv&\hspace*{-1mm} 3\times3&\hspace*{-2mm}\big(\textrm{mod}\;2\big)\\
\multicolumn{4}{l}{\hspace*{-3.2cm}\textrm {\small Rule III,  Table \ref{tab:retrograde-pathway7}\,:}}\\
\hspace*{-1mm} 9+12+2&\hspace*{-1mm}\equiv& \hspace*{-1mm}7\times3&\hspace*{-2mm}\big(\textrm{mod}\;2\big)\\
 \multicolumn{4}{l}{\hspace*{-3.2cm}\textrm {\small Rule II,\,\,\,  Table \ref{tab:retrograde-pathway11}\,:}}\\
\hspace*{-1mm} 9+4+2&\hspace*{-1mm}\equiv& \hspace*{-1mm}7\times11&\hspace*{-2mm}\big(\textrm{mod}\;2\big)
\end{array}\end{equation}  
%on the left, $J_4^{(31)}$ and $J_{12}^{(31)}$  on the right. The domains of the basic operation and the adjoint operation coincide, and the two uses prescribed by condition (\ref{eq:uses})  allow applying both of them.
%\noindent Similar considerations apply to Rule III (Table \ref{tab:retrograde-pathway7}). Even though the domain of $m$ reduces to $\mathbb{N}\,\backslash\{1,2,3,4\}$ for the adjoint operation, the $m$'s   lie high enough, and condition (\ref{eq:uses}) prescribes  two uses,  one applicable to  the basic operation and the other to the adjoint operation. Rule I may appear problematic in that condition (\ref{eq:uses}) prescribes {\em one} use only$-$ the fact that a $G$-tpe or $J$-type semiprime with only 3 as first factor  (Table \ref{tab:retrograde-pathway3}) may nevertheless lead to  semiprime-indexed primes of sizeable iteration  (Eq.\,\ref{eq:iteration}) speaks  in favor of  using the pair's second member. 

\noindent A scoundrel, who thinks evil of it: For order 31, it follows that ${\cal G}^{(31)}/{\cal J}^{(31)}=1/5$ or ${\cal G}^{(31)}/({\cal G}^{(31)}+{\cal J}^{(31)})=1.\bar{6}$, which is less than 1\% off the value published by NASA for the proportion of ordinary matter over the total of matter, ordinary plus dark. 

\medskip {}
\noindent Regarding extrapolation to order 63, one has to keep in mind that the number of semiprime $J_\rho^{(63)}$ with first factor 3, 7 or 11 is not known (at least not until an AI has computed $(J_\rho)^{(63)}\,$)! What is known is that
\begin{equation}\begin{array}{lclcl}
                          G_{4}^{(63]}&=&18633&=&3\times6211,\\ 
                          G_5^{(63)}&=&58791&=&3\times19597    
                                      \;(\;=\;\vert J_{12}^{(31)}\vert), \\
                          G_6^{(63)}&=&189393&=&3\times63131,\\
                          G_{28}^{(63)}&=&519,693,033&=&3\times173,231,011.
                                      \end{array}\label{eq:autonom} \end{equation}
  \noindent The lack of alternatives, in turn, is striking. However,  as  the immersion effect $G_5^{(63)}=\vert J_{12}^{(31)}\vert$ suggests,  there exist  connections between {\em G}-type and {\em J}-type semiprimes, and twelve {\em J\,}-type semiprimes could well turn out pertinent; whether four each with first factor 3, 7 and 11,   leading to  ${\cal G}^{(63)}/{\cal J}^{(63)}=1/5$ or ${\cal G}^{(63)}/({\cal G}^{(63)}+{\cal J}^{(63)})=1.\bar{6}$  even under order 63 conditions, remains an open question. 

\smallskip{}
 \noindent 
%One can only speculate what picture of reality is being conveyed here.
%For instance, that the anterograde pathway  (order 15, Table \ref{tab:anterograde-pathway}) marks an early phase when primordial black holes were exclusively produced.  %Retrograde pathways of order 31,  63 or  higher might then hint at matter that remembers the place of its origin$-$some black hole's event horizon polarizing the immediate fifth  of quasi-stellar surroundings (Tables ) with conventional vs.\,dark charge (Eq.\,\ref{eq:trilinearrational}).
%In that BH phase,  there exists only a boundary layer one unit thick (Rule I)$-$the event horizon itself,  permeable, except for the process of evaporation, only from one side at a time (condition (\ref{eq:uses})).
%In the sequel, more and more matter is accumulated, until at a critical threshold polarization sets in and the black hole becomes a black hole star. In this BHS phase,  a second boundary layer forms  farther outward, possessing, according to condition (\ref{eq:uses}),   two-units thickness  as soon as order 31 is reached and retrograde pathways (Rule II-III) from both sides of that boundary layer can reconnect with the source of charge, normal or ``dark,''   the direction of rotation serving as the discriminant for this difference in charge (Eq.\,\ref{eq:trilinearrational}).
How can the results been interpreted?  Tables \ref{tab:retrograde-pathway3}-\ref{tab:retrograde-pathway11} might describe evolution alongside gravitational reconnection from the perspective of an external observer. At  places marked in blue, information gets frozen in and time stands still. As soon as any gravitational field sheets start to realign, there's a short phase  (in normal black) during which information can flow freely, until it gets reconfined as the sheet finds its new configuration and time slows down again  (next `blue place'). This sequence of events may repeat any number of times until, eventually, a  nuclear way out is reached: in  case of $G$-type paths, in the form of beta decay  (the role of the number 137\,($=\pi(\textcolor{red}{3\times11})$ in beta decay has been discussed in \cite{Merkel1}); in  case of $J$-type paths,
%(like $\textcolor{blue}{3\times19597}$\hspace*{-14.5mm}\color{red}{\ul{\qquad\quad\,\,\,\,\,\,}}\,\color{black})
places  (marked in blue, but underlined in red) where components of atoms disassemble$-$into neutrons, protons and electrons. Gravity wave detection is possibly still in its infancy;  synchronies between  events and patterns as described  above should  be detectable, though, at least  in principle. Table \ref{tab:anterograde-pathway}, on the other hand, shows that there's no intervening $m$ (black)  between $\textcolor{blue}{3\times5}$\hspace*{-7.5mm}\color{red}{\ul{\quad\,\,\,\,\,\,}}\,\,\,\color{black} {and} $\textcolor{red}{3\times7},$ {as maybe merely  describing the evolution of  information from the  perspective of an internal observer, stating that  there's no pause between the  ``bang'' and primordeal nucleosynthesis$-$the opposite end of, or end to, to beta decay.}


\begin{thebibliography}{Cohn}
 \bibitem[Green]{Green}H.S.Green, Phys.Rev. $\mathbf{90}$, 270 (1953)
\bibitem[Merkel]{Merkel}U.Merkel: {\tt arXiv:math/0608423\sl v41}, 2006
\bibitem[Merkel1]{Merkel1}U.Merkel: {\tt arXiv:1309.5272\sl v48}, 2008
\bibitem[Palmer]{Palmer}T.Palmer: {\tt arXiv:quant-ph/2510.02877\sl v3}, 2025

\bibitem[Wang]{Wang}Z.Wang, K.Hazzard: {\tt arXiv:quant-ph/2308.05203\sl v4}, 2025
\bibitem[Veenema]{Veenema}H.G.Veenema et.\,al.\,: {\tt arXiv:2512.11652},2025

\end{thebibliography}
\end{document}